\documentclass{amsproc}

\usepackage{latexsym}
\usepackage{amsmath}
\usepackage{amssymb}
\usepackage{amsfonts}
\usepackage{verbatim}
\usepackage[latin1]{inputenc}
\usepackage{mathrsfs}
\usepackage{amsfonts}
\usepackage{graphicx}
\usepackage{colortbl,dcolumn}
\usepackage{amsmath}
\usepackage{psfrag}
\usepackage{booktabs}
\usepackage{ulem}

\newtheorem{tm}{Theorem}[section]
\newtheorem{rk}{Remark}[section]

\newtheorem{prop}{Proposition}[section]
\newtheorem{lm}{Lemma}[section]
\newtheorem{cor}{Corollary}[section]

\newcommand{\N}{\mathbb N}
\newcommand{\bi}{\mathbf i}
\newcommand{\bs}{\mathbf s}
\newcommand{\E}{\mathbb E}
\newcommand{\R}{\mathbb R}
\newcommand{\PP}{\mathbb P}

\newcommand{\LL}{\mathcal L}

\newcommand{\HH}{\mathbb H}
\newcommand{\FF}{\mathcal F}

\newcommand{\FFF}{\mathscr F}
\newcommand{\OOO} {\mathscr O}

\begin{document}

\title[On  Global Existence and Blow-up for   Damped Stochastic NLS Equation]
{On  Global Existence and Blow-up for   Damped Stochastic Nonlinear Schr\"odinger Equation}

\author{Jianbo Cui}
\address{1. LSEC, ICMSEC, 
Academy of Mathematics and Systems Science, Chinese Academy of Sciences, Beijing,  100190, China\qquad
2. School of Mathematical Science, University of Chinese Academy of Sciences, Beijing, 100049, China}
\curraddr{}
\email{jianbocui@lsec.cc.ac.cn}
\thanks{}
\author{Jialin Hong}
\address{1. LSEC, ICMSEC, 
Academy of Mathematics and Systems Science, Chinese Academy of Sciences, Beijing,  100190, China\qquad
2. School of Mathematical Science, University of Chinese Academy of Sciences, Beijing, 100049, China}
\curraddr{}
\email{hjl@lsec.cc.ac.cn}
\thanks{}
\author{Liying Sun}
\address{1. LSEC, ICMSEC, 
Academy of Mathematics and Systems Science, Chinese Academy of Sciences, Beijing,  100190, China\qquad
2. School of Mathematical Science, University of Chinese Academy of Sciences, Beijing, 100049, China}
\curraddr{}
\email{liyingsun@lsec.cc.ac.cn}
\thanks{}

\subjclass[2010]{Primary 60H35; Secondary 60H15,
60G05}

\keywords{stochastic Schr\"odinger equation,
global existence, 
blow-up,
exponential integrability }

\date{\today}

\dedicatory{}

\begin{abstract}
In this paper, we consider the well-posedness of the weakly damped stochastic nonlinear Schr\"odinger(NLS) equation driven by  multiplicative noise.  First, we show the global existence
 of the unique solution for the damped stochastic NLS equation in critical case. Meanwhile, the exponential integrability  of the solution is proved, which implies the continuous dependence on the initial data. 
Then, we analyze the effect of the damped term and noise on the blow-up phenomenon. By modifying the associated energy, momentum and variance identity, we deduce a sharp blow-up condition for damped stochastic NLS  equation in supercritical case. 
Moreover, we show that when the damped effect is large enough, the damped effect can prevent the  blow-up of the solution with  high probability.
\end{abstract}

\maketitle


\section{Introduction}
The nonlinear Schr\"odinger equation, as one of the basic models for nonlinear waves,  has many physical applications to, e.g. nonlinear optics, plasma physics and quantum  field theory and so on (see e.g. \cite{Bou99,Caz03,BD99,BD03,SS99}). 

In this paper, we consider the weakly damped stochastic NLS equation driven by a linear multiplicative noise
in focusing mass-(super)critical range, 
\begin{align}\label{nls}
du&=\bi(\Delta u+ \lambda |u|^{2\sigma}u) dt -audt+\bi u \circ dW(t), \\\nonumber
u(0)&=u_0, 
\end{align}
where $\frac 2d\le \sigma < \frac 2{(d-2)^+}$, $\lambda =1$, $a\ge0$, $x\in  \R^d$, $t \ge 0$ and $``\circ" $ stands for a Stratonovich product.  Here $W=\{W(t):\ t\in [0,T]\}$ is an $L^2(\mathbb R^d)$-valued $Q$-Wiener process on  a stochastic basis $(\Omega, \FFF, \FFF_t, \PP)$, i.e., there exists an orthonormal basis $\{e_k\}_{k\in \N_+} $ of $L^2(\mathbb R^d)$ and a sequence of mutually independent, real-valued Brownian motions $\{\beta_k\}_{k\in \N_+} $ such that $W(t)=\sum\limits_{k\in \N_+} Q^{\frac12}e_k\beta_k(t)$, $t\in [0,T]$.
We will use frequently the equivalent It\^o form of Eq. \eqref{nls} 
\begin{align*}
du&=\bi(\Delta u+ |u|^{2\sigma}u) dt -(a+\frac 12F_Q)udt+\bi u dW(t), \\\nonumber
u(0)&=u_0
\end{align*}
 with $F_Q:=\sum\limits_{k\in \N^+}(Q^{\frac 12}e_k)^2$.

In recent twenty years, much effort has been devoted to studying the well-posedness  of stochastic NLS equation, see \cite{BRZ16,BM14,BD99,BD03,Hor16} and the references therein.
In \cite{BD99} and \cite{BD03}, the local and global existence of the mild  solution of stochastic NLS equation in $L^2$ and in $\HH^1$ are investigated, respectively.
For stochastic NLS equation on a manifold, \cite{BM14} considers the existence and uniqueness of the solution based on the stochastic Strichartz estimate in UMD Banach space.
\cite{Hor16} shows the local well-posedness in $L^2$ via stochastic Strichartz estimates and the global well-posedness in subcritical case.
With the help of the rescaling transformation, \cite{BRZ16} obtains that the local well-posedness  in $\HH^1$ for $\sigma<\frac 2{(d-2)^+}$ with some decay conditions on the noise, and that the global existence when $\lambda=-1, \sigma<\frac 2{(d-2)^+}$ or  $\lambda=1, \sigma<\frac 2{d}.$
The global existence of the solution in $\HH^2$ is presented in \cite{CHL16b} for one dimensional stochastic NLS equation driven by linear multiplicative noises.
For the global existence of the solution of the NLS equation in critical case, there exist much more results in the deterministic case than those in stochastic case.
For instance,
in the deterministic case, \cite{Wei83} finds a threshold R by the optimal  constant of Gagliardo--Nirenberg's inequality, and proves the global existence of the solution when $\sigma =\frac 2d$ and $\|u_0\|<\|R\|$.
After adding the noise and damped effect,
we wonder whether these effect will influence the existence of a threshold in mass critical case for  the damped stochastic NLS equation, which is one of the main interests in this paper.


The another main interest is to study the damped effect and the noise effect on  the blow-up in the focusing mass-(super)critical case.
It's well known that the solution of deterministic NLS equation with $a=0$ in the focusing mass-(super)critical case will blow up in some finite time when $u_0$ possesses some negative Hamiltonian, see \cite{Bou99,Caz03, SS99} and references therein.
When $a>0$, the damped term has the effect to delay the blow-up, see \cite{OT09,SS99,Tsu84} and references therein. 
For instance, for $\sigma>\frac 2d$,  the blow-up may occur for small values of $a$ (see e.g. \cite{Tsu84}) and large enough values of  $a$ can ensure the global existence of the solution for  $\sigma \ge \frac 2d$ (see e.g. \cite{OT09}). 
In stochastic case, the noise also has an impact on blow-up solutions.
\cite{DD05} shows that the noise effect can accelerate the formation of singularity, and that the solution of Eq. \eqref{nls} with $a=0$ in focusing supercritical case will blow up in a finite time with a positive probability when the variance of  the initial datum is finite.
The blow-up solution for the stochastic NLS equation driven by additive noises is considered in \cite{DD02b}.
When the noise of stochastic NLS equation is non-conservative, \cite{BRZ17}  shows that adding a large multiplicative Gaussian noise
can prevent the blow-up in any finite time with high probability.

Throughout  this paper, we assume that the local-wellposedness of the solution of Eq. \eqref{nls} holds.
The local solution $u(\cdot)$ is defined on a random interval $[0,\tau^*(u_0,\omega))$, where $\tau^*(u_0,\omega)$ is a stopping time such that 
\begin{align*}
\tau^*(u_0,\omega)=+\infty,\quad \text{or} \quad \lim_{t \to \tau^*(u_0,\omega)} \|u(t)\|_{\HH^1}=+\infty.
\end{align*}
First, the evolution of charge and energy of the local solution are introduced.
By using the optimal constant of Gagliardo--Nirenberg's inequality, we show the a priori estimation in $\HH^1$-norm, and prove that the threshold R is unchanged when $a\ge0$, $\sigma=\frac d2$ and initial datum is deterministic. 
Moreover, based on the proved exponential integrability of the solution $u$, i.e.,
\begin{align*}
\sup_{t\in[0,\infty)}\E\Big[\exp\Big(\frac {\|u(t)\|_{\HH^1}^2}{e^{\alpha t}}\Big)\Big]\le C(u_0,a,Q)
\end{align*}
with $\alpha$ depending on $u_0,a$ and $Q,$ we obtain the strong continuous dependence on the initial data in one dimensional case, which is not a trivial property for stochastic partial differential equation with non-global coefficients, see \cite{CHJ13, CHL16b} and references therein. 
We would like to mention that this exponential integrability is useful for  studying  the 
continuous dependence on noises, exponential tail 
estimate of the solution, strong and weak convergence rates of numerical approximations, see \cite{CHJ13,CH17, CHL16b, CHLZ17, HJ14} and references therein. 

Next we consider the influence of damped term and noise on the blow-up.
For the damped stochastic NLS equation, that is, $a>0$, the method used in \cite{DD05} to get the blow-up condition is not available since the variance identity of Eq. \eqref{nls} do not have a polynomial expansion.   
To overcome this difficulty, we modify the energy, momentum and variance identity which is similar to \cite{Tsu84}, and deduce a sharp blow-up condition. Indeed, we show that under some mild assumptions on  $u_0$ and $Q$,  if there exist  $z \ge \frac{4a\sigma}{\sigma d-2} $ and $\bar t$ such that
\begin{align*}
\E\Big[V(u_0)\Big]
+4\bar t\E\Big[G(u_0)\Big]&+(8\bar t^2+\frac 83 z\bar t^3)\E\Big[H(u_0)\Big]\\
&+\Big(\frac 43\bar t^3+\frac 43 z\bar t^4\Big)\E\Big[\|u_0\|^2\Big]\|f_Q\|_{L^{\infty}}\le 0,
\end{align*}
where $f_Q=\sum_{k\in \N^+}|\nabla Q^{\frac 12}e_k|^2$,
then 
\begin{align*}
\PP(\tau^*(u_0)\le \bar t)>0.
\end{align*}
This implies that no matter how large the damped effect is, the blow-up phenomenon  will not disappear. 
We remark that the above blow-up condition can be degenerated to  the blow-up condition in
conservative stochastic case and in the deterministic case.
On the other hand, if the noise satisfies more conditions,   using the rescaling transform idea in \cite{BRZ16},  we prove that when $a\to \infty$ and $\sigma \ge \frac 2d$, for any fixed time T, the blow-up of the solution does not  happen with probability 1.

This paper is organized as follows. In Section 2,  we 
study the evolution of charge and energy, and  show the global existence of the unique solution. In Section 3, 
the modified variance identity is given. Based on it, we obtain a sharp blow-up condition. Furthermore, we prove that when the value of the damped coefficient $a$ becomes large enough, the 
solution does not blow up at any finite time with high probability.  At last, We give a short conclusion in Section \ref{sec;5}.

\section{Global existence of solutions for critical stochastic NLS equations}
\label{sec;2}
In this section, we focus on the global existence and some properties  of  the solution  for Eq. \eqref{nls}.
Throughout this paper, we assume that the local well-posedness for Eq. \eqref{nls} holds.
For the local well-posedness for Eq. \eqref{nls}, we 
refer to \cite{BRZ16,BD03,Hor16} and references therein.
When consider the focusing mass-(super)critical case, \cite{DD02b} proves that the solution of Eq. \eqref{nls} blows up with any initial data for the additive case. 
For the stochastic NLS driven by the multiplicative noise, similar situation happens with any initial datum in 
the super-critical case (see e.g. \cite{BRZ17,DD05}).
This phenomenon is different from the deterministic case, where the solution will blow up in some finite time when $u_0$ possesses some negative Hamiltonian  in the focusing mass-(super)critical case (see e.g. \cite{Bou99, Caz03,SS99}).
However, it is still not clear on whether or not
the solution of  Eq. \eqref{nls} equation globally exists in critical case.
Notice that when the noise is independent of space and $a=0$, i.e.,
\begin{align*}
du&=\bi(\Delta u+ |u|^{2\sigma}u) dt - \frac 12udt +\bi ud\beta(t), \\\nonumber
u(0)&=u_0,
\end{align*}
the global existence and blow-up results become more clear.
In this case, one can use the  infinite dimensional Doss--Sussman type transformation $u(t)
=\exp(\bi \beta(t))y(t)$ to get the well-posedness and blow-up results, where $y(t)$ satisfies 
\begin{align*}
dy&=\bi \Delta ydt +\bi |y|^{2\sigma}ydt,\\
y(0)&=u_0.
\end{align*}

We first study the global exsitence of the solution of 
Eq. \eqref{nls} in  the focusing critical case. 
This suggests that  the critical nonlinearity in multiplicative cases is different from the supercritical nonlinearity, and that the critical nonlinearity combined with the dispersion term dominates the behavior of the solution.
For convenience, we assume that  $u_0\in \HH^1$ is a deterministic function and that
$\sum_{k}\|Q^{\frac 12}e_k\|_{\HH^1}^2+
\|f_Q\|_{L^{\infty}}<\infty$ 
with $f_Q=\sum_{k}|\nabla Q^{\frac 12}e_k|^2$. 
The solution $u(\cdot)$ of Eq. \eqref{nls} is defined on a random interval $[0,\tau^*(u_0,\omega))$, where $\tau^*(u_0,\omega)$ is a stopping time such that 
\begin{align*}
\tau^*(u_0,\omega)=+\infty,\quad \text{or} \quad \lim_{t \to \tau^*(u_0,\omega)} \|u(t)\|_{\HH^1}=+\infty.
\end{align*}
To get a priori estimate of $u$, we first study the 
 evolution of 
charge $M(u(t)):=\|u(t)\|^2$ and energy $H(u):=\frac 12\|\nabla u\|^2-\frac \lambda {2\sigma+2}\|u\|_{L^{2\sigma+2}}^{2\sigma+2}$ in the following lemma.

\begin{lm}\label{char}
Assume that $u_0 \in \HH^1$  and $\sum_{k}\|Q^{\frac 12}e_k\|_{\HH^1}^2+
\|f_Q\|_{L^{\infty}
}<\infty$.
For any  $\tau<\tau^*(u_0)$, we 
have 
\begin{align}\label{cha}
M(u(\tau))=e^{-2a\tau}M(u_0), \quad \text{a.s.,}
\end{align}
and  
\begin{align}\label{ham}
H(u(\tau))&=H(u_0)-a\int_{0}^{\tau}(\|\nabla u(s)\|^2-\|u(s)\|^{2\sigma+2}_{L^{2\sigma+2}} )ds\\\nonumber
&\quad-\text{\rm Im}\int_{\R^d}\int_0^{\tau}\bar u(s)\nabla u(s)\nabla dW(s)dx\\\nonumber
&\quad+\frac12 \sum_{k\in \N^+}\int_0^{\tau} \|u(s)\nabla Q^{\frac 12}e_k\|^2ds, \quad \text{a.s.}
\end{align}

\end{lm}
\textbf{Proof}
 For purpose of obtaining the charge and energy evolution of $u$, the truncated argument in \cite{BD03} is applied. In detail, 
let $N\in \N^+$ and $K>0$ and define the operators $\Theta_N,  N \in \N$ by 
\begin{align*}
\Theta_Nv:=\FF^{-1}\Big(\theta(\frac {|\cdot|}{j})*N \Big),
\end{align*}
where $\FF$ is the Fourier transform 
and $\theta \in C_c^{\infty}$ is a real-valued and nonnegative function satisfying $\theta(x) = 1$ for $|x| \le 1$, $\theta(x) = 0$ for $|x|>2$.
Using the above notation, we have the truncated approximation, for $m=(m_1,m_2)\in \N^2$,
\begin{equation}
\begin{split}
\label{spe}
du_K^m&=\bi
\left(\Theta_{m_1} \Delta u_K^m+ \theta\Big{(}\frac {\|u_K^m\|_{\HH^1}}{K}\Big{)} \Theta_{m_2}(|u_K^m|^{2\sigma}u_K^m)- au^m_K- \frac {F_{Q_{m_2}}} {2} u^m_K\right) dt\\
&\quad+\bi u_K^m\Theta_{m_2}dW(t),
\end{split}
\end{equation}
where $F_{Q_{m_2}}:=\sum_{k\in\N^+}(\Theta_{m_2}(Q^{\frac 12})e_k)^2$.
Combining with It\^o formula in $[0,\tau]$ and taking limits as $m \to \infty $, the evolution of the charge \eqref{cha} is obtained by choosing a large enough $K$.
Similarly, using the above arguments, the energy evolution law \eqref{ham} can be proved. 
\qed

\begin{rk}
The truncated argument is also available for stochastic 
NLS equation with  the homogenous Dirichlet boundary condition. 
In this case,  replacing $\Theta_N$ by the projection operator $P^N$,
then the truncated Galerkin approximated equation becomes
\begin{equation}
\begin{split}
\label{spe}
du_K^N&=\bi\left(\Delta u_K^N
+ \theta\Big(\frac {\|u_K^N\|_{\HH^1}}{K}\Big) P^N(|u_K^N|^{2\sigma}u_K^N)- aP^Nu^N_K
- P^N\Big(\frac {F_Q} {2} u^N_K\Big)\right) dt\\
&\quad  +\bi P^N(u_K^NdW(t)),
\end{split}
\end{equation}
where $K>0, N\in \N^+$. 
The inverse inequality, for $\bs \ge 1$,
\begin{align*}
\|u_K^N\|_{\HH^{\bs}}
\le C(N) \|u_K^N\|_{\HH^1}
\end{align*}
implies the coefficients of Eq. \eqref{spe} are globally Lipschitz.
Therefore, by the arguments in \cite{BD03},  the result of Lemma \ref{char} holds. 
\end{rk}

In order to illustrate the global well-posedness result, we introduce the optimal constant for Gagliardo--Nirenberg inequality  and its corresponding ground state solution
(see e.g. \cite{Wei83}).

\begin{lm}\label{gne}
The best constant $C_{\sigma, d}$ for  Gagliardo--Nirenberg inequality 
\begin{align}\label{gn}
\|f\|_{L^{2\sigma+2}}^{2\sigma+2}\le C_{\sigma, d}
\|\nabla f\|^{\sigma d}\|f\|^{2+\sigma(2-d)}
\end{align}
with $f\in \HH^1(\R^d) $, $0<\sigma<\frac 2{d-2}$ and $d\ge 2$ is given by 
\begin{align*}
C_{\sigma, d}=(\sigma+1)\frac {2(2+2\sigma-\sigma d)^{-1+\frac {\sigma d}2}}{(\sigma d)^{\frac {\sigma d}2}}\frac 1{\|R\|^{2\sigma}},
\end{align*}
where $R$ is the positive solution (ground state solution) 
of $\Delta R-R+R^{2\sigma+1}=0$.
\end{lm}

Based on Lemma \ref{char}  and Lemma \ref{gne}, we are in position to show the global existence of $u$.
For the sake of simplicity, 
the procedure about truncated arguments and taking limits is omitted in the rest of this paper.

\begin{tm}\label{well}
 Assume that $u_0\in \HH^1$ with $\|u_0\|< \|R\|$,
 $\sum_{k}\|Q^{\frac 12}e_k\|_{\HH^1}^2+$ $\|f_Q\|_{L^{\infty}}<\infty$. Then there exists a unique global 
solution of Eq. \eqref{nls} in $\HH^1$, i.e., $\tau^*(u_0)=\infty$.
\end{tm}
\textbf{Proof}
Since the local well-posedness of Eq. \eqref{nls} is shown  by \cite{BD03, BM14,BRZ16}, we only need to get the uniform boundedness of $\|u\|_{\HH^1}$ 
to ensure the global existence of the solution.
By the charge conservation law, Gagliardo--Nirenberg inequality and $\sigma d=2,$
\begin{align*}
\|\nabla u\|^2\le 2H(u)+\frac 1{\sigma+1}C_{\sigma, d}\|u\|^{2\sigma}\|\nabla u\|^2,
\end{align*}
where $C_{\sigma, d}=\frac {\sigma+1}{\|R\|^{2\sigma}}$ and $R$ is the ground state solution of $\Delta R-R +R^{2\sigma+1}=0$.
Then the energy evolution of $u$ implies that for any $T_0>0$, any stopping time $\tau < \inf(T_0,\tau^*(u_0))$ and any time $t\le \tau$,
\begin{align*}
\left(1-\frac {\|u(t)\|^{2\sigma}}{\|R\|^{2\sigma}}\right)\|\nabla u(t)\|^2
& \le 2H(u_0)-2a\int_{0}^{t}(\|\nabla u(s)\|^2-\|u(s)\|^{2\sigma+2}_{L^{2\sigma+2}} )ds
\\
&\quad-2\text{Im}\int_{\OOO}\int_0^{t}\bar u(s)\nabla u(s)\nabla dW(s)dx\\
&\quad+\sum_{k\in \N^+}\int_0^{t}\int_{\OOO}|u(s)|^2|\nabla Q^{\frac 12}e_k|^2dxds,
\end{align*}
where  $H(u_0)> 0$ since  $\|u_0\|<\|R\|$.
After taking expectation, we have 
\begin{align*}
\left(1-\frac {\|u_0\|^{2\sigma}}{\|R\|^{2\sigma}}\right)
\E\Big[\|\nabla u(t)\|^2\Big]
&\le 2H(u_0)-4a\int_0^t\E \Big[H(u(s))\Big]ds\\
&\quad+\int_0^t \frac {2a\sigma}{\sigma+1}\E \Big[\|u(s)\|^{2\sigma+2}_{L^{2\sigma+2}}\Big]ds\\
&\quad+ \E \left(\sum_{k\in \N^+}\int_0^{t}\int_{\OOO}|u(s)|^2|\nabla Q^{\frac 12}e_k|^2dxds\right).
\end{align*}
{Then  H\"older inequality and Sobolev embedding theorem yield that}
\begin{align*}
\quad&\left(1-\frac{\|u_0\|^{2\sigma}}{\|R\|^{2\sigma}}\right)
\E\Big[\|\nabla u(t)\|^2\Big]\\
&\le 2H(u_0)-4a\int_0^t\E \Big[\|\nabla u(s)\|^2\Big]ds+ \frac {2a\sigma}{\sigma+1} C_{\sigma,d}\int_0^t \E \Big[ \|\nabla u(s)\|^2\|u(s)\|^{2\sigma}\Big]ds \\
&\quad+ \int_0^t\sum_{k\in \N^+}\|\nabla Q^{\frac 12}e_k\|_{L^{\infty}}^2\E\Big[\|u(s)\|^2\Big]ds\\
&\le 2H(u_0)+
\frac {2a\sigma}{\sigma+1}C_{\sigma,d}\int_0^t e^{-2as}\E\Big[\|\nabla u(s)\|^2\Big]\|u_0\|^{2\sigma}ds
\\
&\quad+\int_0^{t}\sum_{k\in \N^+}\|\nabla Q^{\frac 12}e_k\|_{L^{\infty}}^2e^{-2as}\|u_0\|^2ds.
\end{align*}
Gronwall inequality implies that 
\begin{align*}
\sup_{t\le \tau}\E\Big[\|\nabla u(t)\|^2\Big]\le C(u_0,R,a,Q,T_0).
\end{align*}
Moreover, using Burkholder--Davis--Gundy inequality and Young inequality,
\begin{align*}
\quad&\E\Big[\sup_{t\le \tau}\|\nabla u(t)\|^2\Big]\\
&\le C(u_0,R,a,Q,T_0)+C\E\Big[\sup_{t\le \tau}\Big|\int_0^t \int_{\R^d} \bar u(s)\nabla u(s)\nabla dW(t)dx\Big|\Big]\\
&\le C(u_0,R,a,Q,T_0)+C\E\Big[\Big|\int_0^\tau \sum_{k\in \N^+}\| u(s)\nabla u(s)\nabla Q^{\frac 12}e_k\|^2dt\Big|^{\frac 12}\Big]\\
& \le C(u_0,R,a,Q,T_0)+\epsilon\E\Big[\sup_{t\le \tau}\|\nabla u(t)\|^2\Big]\\
\quad\quad&\quad+C(\epsilon) \E\Big[ \int_0^{\tau}\sum_{k\in \N^+}\|\nabla Q^{\frac 12}e_k\|_{L^{\infty}}^2 \|u(s)\|^2dt\Big],
\end{align*} 
where $\epsilon<\frac 12$. The last term of the above inequality is bounded by H\"older inequality and the charge evolution law \eqref{cha}, which in turns implies that
$$\E\Big[\sup_{t\le \tau}\|\nabla u(t)\|^2\Big]\le C(u_0,R,a,Q,T_0).$$ 
These a priori estimations combined with local well-posedness imply the global existence of the unique solution. 
\qed
\\

The  condition $\|u_0\|< \|R\|$  in Theorem \ref{well} is a sufficient condition for the global existence of the solution. Based on it, we get a upper bound of the probability of the  blow-up with a random initial datum at any finite time.

\begin{cor}\label{nes}
Assume that $u_0$ is random initial datum. Under the condition of  Theorem \ref{well}, we have 
\begin{align*}
\PP(\tau^*(u_0)<\infty)\le \PP(\|u_0\|\ge\|R\|).
\end{align*} 
\end{cor}

Moreover, we can get the following exponential integrability of $u$, which is useful for studying the 
continuous dependence on initial data and noise, exponential tail 
estimate of the solution, strong and weak convergence rates of numerical approximations (see e.g. \cite{ CH17,CHL16b, CHLZ17}).
  
\begin{prop}\label{exp}
Assume that 
the condition of Theorem \ref{well} holds, then there exist  constants $C=C(u_0,R,a,Q) $ and $\alpha=\alpha(u_0,R,a,Q)$  such that 
\begin{align*} 
\sup_{t\in [0,T]}\E\Big[\exp\left(\frac {\|u\|^2_{\HH^1}}{e^{\alpha t}}\right)\Big]\le C.
\end{align*}
\end{prop}
\textbf{Proof}
Denote by a positive number $c:=1-\frac {\|u_0\|^{2\sigma}}{\|R\|^{2\sigma}}<1$,
$\mu(u)=\bi \Delta u +\bi |u|^2u-\frac 12F_Qu-a$ and 
$\sigma(u)=\bi uQ^{\frac 12}$.
Based on Gagliardo--Nirenberg inequality \eqref{gn}, we only need to show the boundedness of $\sup\limits_{t\in[0,T]}\E\Big[\exp\Big(\frac {2H(u(t))}{ce^{\alpha t}}\Big)\Big]$.
 Using the truncated arguments and taking limits, we get
\begin{align*}
&\frac {2 DH(u(s))\mu(u(s))}{c}+ \frac {\text{tr}[D^2H(u(s))\sigma(u(s))\sigma^*(u(s))]}c
+ \frac {4\|\sigma^*(u(s))DH(u(s))\|^2}{c^2e^{\alpha t}}\\
&=\frac 1c{\sum\limits_k\|u(s)\nabla Q^{\frac 12}e_k\|^2}
-\frac {2a}c\|\nabla u(s)\|^2+\frac {2a}c\|u(s)\|^{2\sigma+2}_{L^{2\sigma+2}}\\
&\quad+\frac {4}{c^2e^{\alpha t}}
\sum_k\|\nabla u(s)\|^2\|u(s)\nabla Q^{\frac 12}e_k\|^2.
\end{align*}
Again by the Gagliardo--Nirenberg inequality \eqref{gn} and  the Sobolev embedding theorem,  we obtain
\begin{align*}
&\frac {2 DH(u(s))\mu(u(s))}{c}+ \frac {\text{tr}[D^2H(u(s))\sigma(u(s))\sigma^*(u(s))]}c
+ \frac {4\|\sigma^*(u(s))DH(u(s))\|^2}{c^2e^{\alpha t}}\\
&\le 
\frac 1c\|u(s)\|^2\|f_Q\|_{L^{\infty}}
-\frac {4a}cH(u(s))+\frac {2a\sigma}{c(\sigma+1)}C_{\sigma, d}\|\nabla u(s)\|^2\|u(s)\|^{2\sigma}\\
&\quad+\frac {4}{c^2e^{\alpha t}}\|\nabla u(s)\|^2\|u(s)\|^2
\|f_Q\|_{L^{\infty}}\\
&\le \frac 1c\|u(s)\|^2\|f_Q\|_{L^{\infty}}
-\frac {4a}cH(u(s))+\frac {4a\sigma}{c^2(\sigma+1)}\|u(s)\|^{2\sigma}C_{\sigma, d}H(u(s))\\
&\quad+\frac {8}{c^3e^{\alpha t}}\|u(s)\|^2
\|f_Q\|_{L^{\infty}}H(u(s)).
\end{align*}
Charge evolution law in Lemma \ref{char} leads that
\begin{align*}
&\frac {2 DH(u(s))\mu(u(s))}{c}+ \frac {\text{tr}[D^2H(u(s))\sigma(u)\sigma^*(u)]}c
+ \frac {4\|\sigma^*(u)DH(u(s))\|^2}{c^2e^{\alpha t}}\\
&\le \frac {e^{-2at}}c\|u_0\|^2\|f_Q\|_{L^{\infty}}\\
&\quad+\Bigg(-2a+\frac {2a\sigma}{c(\sigma+1)}\|u_0\|^{2\sigma}C_{\sigma, d}
+\frac {4}{c^2e^{(\alpha+ 2a)t}}\|u_0\|^2
\|f_Q\|_{L^{\infty}}\Bigg)\frac 2cH(u (s)).
\end{align*}
Then Lemma 3.1 in \cite{CHL16b} implies that 
\begin{align*}
\sup_{t\in [0,T]} \E\Bigg[\exp\Bigg(\frac {\frac 2cH(u(t))}
{e^{\alpha t}}\Bigg)\Bigg]\le \E\Big[\exp\Big(\frac 2cH(u_0)\Big)\Big]\exp\Big(\int_0^T\frac {\|u_0\|^2\|f_Q\|_{L^{\infty}}}{e^{(\alpha+2a) t}}dt\Big),
\end{align*}
where $\alpha\ge -2a+\frac {2a\sigma}{c(\sigma+1)}\|u_0\|^{2\sigma}C_{\sigma, d}+\frac {4}{c^2}\|u_0\|^2
\|f_Q\|_{L^{\infty}}^2$. 
\qed
\\

Applying the above exponential integrability of exact solution,
we deduce the following strongly continuous dependence on the initial data.

\begin{cor}
 Assume that $d=1, \sigma=2$, $a>0$, $u_0,v_0 \in \HH^1$ and $\max(\|u_0\|,\|v_0\|)<\min \Big(\frac {4ac(2c-1)}{4c^2+\|f_Q\|_{L^{\infty}}}, \frac {ac(2c-1)}{\|f_Q\|_{\infty}}\Big)$ with $\frac 12<c<1$. Let
$u=\{u(t): t\in[0,T]\}$ and $v=\{v(t): t\in [0,T]\}$ be the solutions of Eq.\eqref{nls} with initial data $u_0$ and $v_0 $, respectively. Under the condition of Proposition \ref{exp},
then there exists a constant $C = C(p,u_0,v_0,Q,a)$ such that
\begin{align*}
\sup_{t\in [0,T]}\E\Big[\|u(t)-v(t)\|^2\Big]
\le C\E\Big[\|u_0-v_0\|^2\Big].
\end{align*}
\end{cor}
\textbf{Proof}
Applying the truncated arguments, It\^o fromula  and taking limits yield that
\begin{align*}
\|u(t)-v(t)\|^2&=\|u_0-v_0\|^2-2a\int_0^t\|u(s)-v(s)\|^2ds\\
&\quad+2\int_0^t\<u(s)-v(s),\bi\Delta(u(s)-v(s))\>ds 
\\
&\quad+2\int_0^t\<u(s)-v(s),\bi(|u(s)|^4u(s)-|v(s)|^4v(s))\>ds\\
&=\|u_0-v_0\|^2-2a\int_0^t\|u(s)-v(s)\|^2ds\\
&\quad+2\int_0^t\<u(s)-v(s),\bi(|u(s)|^4u(s)-|v(s)|^4v(s))\>ds,
\end{align*}
Since for $a,b \in \mathbb C$, $|a|^4a-|b|^4b=(|a|^4+|b|^4)(a-b)+ab(|a|^2+|b|^2)(\bar a-\bar b)+|a|^2|b|^2(a-b)$, combining with Young inequality and Gagliardo--Nirenberg inequality, we get
\begin{align*}
\|u(t)-v(t)\|^2&=\|u_0-v_0\|^2-2a\int_0^t\|u(s)-v(s)\|^2ds\\
&\quad+2\int_0^t\<u(s)-v(s),\bi u(s)v(s)(|u(s)|^2+|v(s)|^2)(\bar u(s)-\bar v(s))\>ds\\
&\le \|u_0-v_0\|^2-2a\int_0^t\|u(s)-v(s)\|^2ds\\
&\quad+2\int_0^t\Big(\|u(s)\|^4_{L^{\infty}}+\|v(s)\|^4_{L^{\infty}}\Big)\|u(s)-v(s)\|^2ds.
\end{align*}
Gronwall inequality and  Gagliardo--Nirenberg inequality lead that
\begin{align*}
\|u(t)-v(t)\|^2
&\le \exp\Big(\int_0^t 2(\|u(s)\|^4_{L^{\infty}}+\|v(s)\|^4_{L^{\infty}}) ds\Big)\|u_0-v_0\|^2\\
&\le \exp\Big(\int_0^t 8e^{-2as}\|u_0\|^2\|\nabla u(s)\|^2+8e^{-2as}\|v_0\|^2\|\nabla v(s)\|^2 ds\Big)\|u_0-v_0\|^2.
\end{align*}
After taking expectation, by Young inequality, we obtain
\begin{align}\label{cont}
\E\Big[\sup_{t\in [0,T]}\|u(t)-v(t)\|^2\Big]
&\le \sqrt{\E\Big[\exp\Big(\int_0^T 16e^{-2as}\|u_0\|^2\|\nabla u(s)\|^2ds\Big)\Big]}\\\nonumber 
&\quad\times \sqrt{\E\Big[\exp\Big(\int_0^T 16e^{-2as}\|v_0\|^2\|\nabla v(s)\|^2 ds\Big)\Big]}
\|u_0-v_0\|^2
\end{align}
 It is obvious that if these two exponential moments in the above inequality are bounded, the theorem is proved.
 For simplicity, we  take $\E\Big[\exp\Big(\int_0^T 16e^{-2as}\|u_0\|^2\|\nabla u(s)\|^2ds\Big)\Big]$ as example.
By Jensen inequality, 
for $\alpha<2a$,
\begin{align*}
&\E\Big[\exp\Big(\int_0^T 16e^{-2as}\|u_0\|^2\|\nabla u(s)\|^2ds\Big)\Big]\\
&=\E\Big[\exp\Big(\int_0^T e^{-(2a-\alpha)s}16\|u_0\|^2e^{-\alpha s}\|\nabla u(s)\|^2ds\Big)\Big]\\
&\le \E\Big[\exp\Big( \frac {16}{2a-\alpha}\|u_0\|^2e^{-\alpha s}\|\nabla u(s)\|^2\Big)\Big].
\end{align*}
Then we take $\alpha=-2a+\frac {2a\sigma}{c(\sigma+1)}\|u_0\|^{2\sigma}C_{\sigma, d}+\frac {4}{c^2}\|u_0\|^2
\|f_Q\|_{L^{\infty}}$
such that 
$\frac {16}{2a-\alpha}\|u_0\|^2\le 1$ and $\alpha<2a$. 
Indeed, the assumption on $u_0$ and $v_0$, together with Gagliardo--Nirenberg inequality Eq. \eqref{gn}, implies
\begin{align*}
\|u_0\|^2&\le
\frac a2-\frac {a}{4c}-\frac {1}{4c^2}\|u_0\|^2
\|f_Q\|_{\infty}\\
&=\frac {a(2c-1)}{4c}-\frac {1}{4c^2}\|u_0\|^2
\|f_Q\|_{\infty}. 
\end{align*}
and
\begin{align*}
\frac 1c+\frac 1{ac^2}\|u_0\|^2
\|f_Q\|_{\infty}< 2,
\end{align*}
which ensure that $\frac {16}{2a-\alpha}\|u_0\|^2\le 1$ and $\alpha<2a$.
Proposition \ref{exp}, combining the above estimations, implies the uniform boundedness of the exponential moments of Eq. \eqref{cont}, which 
completes the proof. 
\qed
\begin{rk}
Due to the fact that the ground state solution R(x) in one dimension is $\frac {3^{\frac 14}}{\sqrt{\cosh(2x)}}$ with $\|R\|^2=\pi\sqrt{3}$, 
the above strong continuous dependence result on initial data holds with $\max (\|u_0\|,\|v_0\|)< \sqrt[4]{\frac {3\pi^2}2}$ when $a$ becomes large enough.
For $ \|u_0\|< \|R\|$, $a\ge 0$, by the above arguments, we can get the continuous dependence on initial data in pathwise sense,
\begin{align*}
\sup_{t\in[0,T]}\|u(t)-v(t)\|
\le C(\omega)\|u_0-v_0\|.
\end{align*}
\end{rk}

\section{Blow-up of solutions in focusing mass-(super)critical case}
\label{sec;3}

As shown in Section \ref{sec;2}, the result about well-posedness for Eq. \eqref{nls} in the critical case is similar to that in deterministic case. 
Notice that this phenomenon is  different from  the additive case with $\sigma \ge \frac d2$ and the multiplicative case with $\sigma >\frac d2$ (see e.g. \cite{DD02b,DD05}), where the singularity happens 
in any finite time with a positive probability for any initial datum.

In fact, the authors in \cite{DD05} show that 
for $\sigma\ge \frac 2d$,  if  $u_0\in L^2(\Omega; \Sigma)\cap L^{2\sigma+2}(\Omega;L^{2\sigma+2}(\R^d))$, $f_Q=\sum_{k\in\N^+}|\nabla Q^{\frac 12}e_k|^2$ and
for some $\bar t>0$,
\begin{align*}
\E\Big[V(u_0)\Big]+4\E\Big[G(u_0)\Big]\bar t+8\E\Big[H(u_0)\Big]\bar t^2+\frac 43\bar t^3\|f_Q\|_{L^{\infty}}\E\Big[M(u_0)\Big]<0,
\end{align*}
then $\PP(\tau^{*}(u_0)\le\bar t)>0$.
The above result implies that if the energy of $u_0$ is a.s. negative, then  $\PP(\tau^{*}(u_0)\le t)>0$ for some $t>0$ provided the noise is not too strong, i.e., $\|f_Q\|_{L^{\infty}}$ is small enough.
The natural question is whether the damped effect can prevent the blow-up phenomenon or not in stochastic case. 

To study the  blow-up phenomenon, we introduce the 
finite variance space
\begin{align*}
\Sigma=\{v\in \HH^{1}: |x|v\in \HH \}
\end{align*}
endowed with the norm $\|\cdot\|_{\Sigma}$:
\begin{align*}
\|v\|_{\Sigma}^2=\||x|v\|^2+\|v\|_{\HH^1}^2,
\end{align*}
the variance 
\begin{align*}
V(v)=\int_{\R^d}|x|^2|v(x)|^2dx,\quad  v\in \Sigma
\end{align*}
and the momentum
\begin{align*}
G(v)=\text{Im} \int_{\R^d} \bar v(x) x\cdot \nabla v(x)dx, 
\quad  v\in \Sigma.
\end{align*}
With the help of a smoothing procedure and truncated arguments (see e.g. \cite{DD05}),  we can prove rigorously the evolution laws of $V$ and $G$ for the damped stochastic NLS equation.

\begin{prop}\label{vg}
Assume that $u_0\in \Sigma$. Under the conditions of Lemma \ref{char},   for any stopping time $\tau<\tau^*(u_0)$ a.s., we have 
\begin{align*}
V(u(\tau))=V(u_0)+4\int_0^\tau G(u(s))ds-2a\int_{0}^\tau V(u(s))ds,
\end{align*}
and
\begin{align*}
G(u(\tau))&= G(u_0)+4\int_{0}^{\tau}H(u(s))ds-2a\int_0^\tau G(u(s))ds+\frac {2-\sigma d}{\sigma+1}\int_0^\tau \|u(s)\|_{L^{2\sigma+2}}^{2\sigma+2}ds\\
&\quad+\sum_{k\in \N}\int_0^\tau\int_{\R^d}|u(s,x)|^2x\cdot \nabla (Q^{\frac 12}e_k)(x)dxd\beta_k(s).
\end{align*}

\end{prop}
\textbf{Proof}
Applying It\^o formula to $V$ and $G$, integration by parts and taking the imaginary part of the integration, we obtain 
\begin{align*}
V(u(\tau))=V(u_0)+4\int_{0}^\tau G(u(s))ds-2a\int_{0}^\tau V(u(s))ds,
\end{align*}
and
\begin{align*}
G(u(\tau))&=G(u_0)
+2\int_{0}^\tau \text{Im}\int_{\R^d} x\cdot
\nabla u\Big(-\bi\Delta \bar u -\bi |u|^{2\sigma}\bar u-a\bar u-\frac12 F_Q\bar u\Big)dxds\\ 
&\quad-d\int_{0}^\tau \text{Im}\int_{\R^d} 
\Big(\bi\Delta u +\bi |u|^{2\sigma} u-au-\frac12 F_Q u\Big)\bar u dx ds\\
&\quad+2\int_{0}^\tau \text{Im}\int_{\R^d}-\bi x\nabla u \bar udxdW(s)-d\int_{0}^\tau\text{Im}\int_{\R^d}\bi |u|^2dxdW(s)\\
&\quad+\int_{0}^\tau\text{Im}\int_{\R^d}x\cdot \nabla(uQ^{\frac 12}e_k)Q^{\frac 12}e_k\bar udxds\\
&=G(u_0)+d\int_0^\tau\Big(\|\nabla u(s)\|^2-\|u(s)\|^{2\sigma+2}_{L^{2\sigma+2}}\Big)ds\\
&\quad-2\int_0^\tau\text{Im}\int_{\R^d}x\cdot
\nabla u\bi\Delta \bar udxds\\
&\quad-2\int_0^\tau\text{Im}\int_{\R^d}x\cdot \nabla u (\bi |u|^{2\sigma+2}\bar u)dxds-2a\int_0^\tau G(u(s))ds\\
&\quad+\sum_{k\in \N}\int_0^\tau\int_{\R^d}|u(s,x)|^2x\cdot \nabla (Q^{\frac 12}e_k)(x)dxd\beta_k(s).
\end{align*}
By the definition of $H$ and $\sigma d=2$, we get 
\begin{align*}
G(u_\tau)&= G(u_0)+4\int_0^\tau H(u(s))ds
+\frac {2-\sigma d}{\sigma+1}\int_0^\tau \|u(s)\|_{L^{2\sigma+2}}^{2\sigma+2}ds-2a\int_0^\tau G(u(s))ds\\
&\quad+\sum_{k\in \N}\int_0^\tau\int_{\R^d}|u(s,x)|^2x\cdot \nabla (Q^{\frac 12}e_k)(x)dxd\beta_k(s).
\end{align*}
\qed

For the damped stochastic NLS equation, the method in \cite{DD05} is not available
since the damped effect will lead that the expansion of 
$V$ produces  many addition terms which can not be estimated directly. We introduce the modified energy, 
invariance and momentum as in \cite{Tsu84} and study the evolution of these modified quantities to investigate  the  blow-up condition for supercritical case $\sigma d>2, a>0$.

\begin{lm}\label{mod}
Let $b\in \R$. Under the same condition of Proposition \ref{vg}, for any stopping time $\tau< \tau^*(u_0)$,
we have 
\begin{align}\nonumber
e^{b\tau}H(u(\tau))&=H(u_0)+b\int_{0}^{\tau}e^{bs}H(u(s))ds-a\int_{0}^{\tau}e^{bs}\Big(\|\nabla u(s)\|^2-\|u(s)\|^{2\sigma+2}_{L^{2\sigma+2}} \Big)ds\\\label{bh}
&\quad-\text{\rm Im}\int_{\OOO}\int_0^{\tau}e^{bs}\bar u\nabla u\nabla dWdx
+\frac12 \sum_{k\in \N^+}\int_0^{\tau} e^{bs}\|u\nabla Q^{\frac 12}e_k\|^2ds, \\\nonumber
e^{b\tau}G(u(\tau))&= G(u_0)-(2a-b)\int_0^\tau e^{bs}G(u(s))ds\\\label{bg}
&\quad+2\int_{0}^{\tau}e^{bs}\Big(2H(u(s))+\frac {2-\sigma d}{2\sigma+2}\|u(s)\|^{2\sigma+2}_{L^{2\sigma+2}}\Big)ds\\\nonumber 
&\quad+ \sum_{k\in \N}\int_0^\tau\int_{\R^d}e^{bs}|u(s,x)|^2x\cdot \nabla (Q^{\frac 12}e_k)(x)dxd\beta_k(s),
\end{align}
and
\begin{align}\label{bv}
e^{b\tau}V(u(\tau))&=V(u_0)+4\int_0^\tau e^{bs}G(u(s))ds-(2a-b)\int_{0}^\tau e^{bs}V(u(s))ds.
\end{align}
\end{lm}
\textbf{Proof}
The proof is similar to  the proof of Lemma \ref{char}
and Proposition \ref{vg} by using  smoothing procedures, truncated arguments, integration by parts and It\^o formula.  More details, we refer to \cite{DD05}.
\qed
\\

Based on Lemma \ref{mod}, we prove a preliminary result on the blow-up condition for Eq. \eqref{nls} in the supercritical case.

\begin{prop}\label{blow}
Let  $\sigma d>2,$ $a>0$, $b<2a$ satisfy $\frac {4a\sigma}{\sigma d-2}\le2a-b,$  $u_0\in L^2(\Omega; \Sigma)\cap L^{2\sigma+2}(\Omega; L^{2\sigma+2})$ and $\sum_{k\in \N^+}\|Q^{\frac 12}e_k\|_{\HH^1}^2+\|f_Q\|_{L^{\infty}}<\infty$.
Assume in addition that for some $0<y \le \frac 1{2a-b}$ such that
\begin{align}\label{bcon}
\E\Big[V(u_0)+4y G(u_0)+16y^2H(u_0)+
8y^3\|u_0\|^2\|f_Q\|_{\infty}\Big]< 0,
\end{align}
 then for some $\bar t$, we have 
\begin{align*}
\PP(\tau^*(u_0)\le \bar t)>0.
\end{align*}
\end{prop}
\textbf{Proof}
We prove the assertion by contradiction. Assume that the solution $u$ exists globally.
Then for any $t>0$, $\tau^*(u_0)> t$ a.s. Then we take 
$\tau=t$.  The evolution law of modified energy Eq. \eqref{bh} , charge evolution law Eq. \eqref{cha} and taking expectation leads that 
\begin{align*}
&\quad\E \Big[e^{bt}H(u( t))\Big]\\
&=H(u_0)+b\int_{0}^{ t}\E\Big[e^{bs}H(u(s))\Big]ds-a\int_{0}^{ t}\E\Big[e^{bs}\Big(\|\nabla u(s)\|^2-\|u(s)\|^{2\sigma+2}_{L^{2\sigma+2}} \Big)\Big]ds\\
&\quad+\frac12 \sum_{k\in \N^+}\int_0^{ t} \E\Big[e^{bs}\|u\nabla Q^{\frac 12}e_k\|^2\Big]ds\\
&\le H(u_0)+
(\frac b2-a)\int_{0}^{t}e^{bs}\E\Big[ \|\nabla u\|^2
+\frac 2{b-2a}(a-\frac b{2\sigma+2})\|u\|_{L^{2\sigma+2}}^{2\sigma+2} \Big]ds\\
&\quad+\frac12 \int_{0}^{t}\E\Big[e^{(b-2a)s}\|u_0\|^2\|f_Q\|_{L^{\infty}}\Big]ds.
\end{align*}
Next, we aim to show a priori estimate on $e^{bt}G(u(t))$. For simplicity, we denote $\tilde H(u):=\|\nabla u\|^2-\frac {\sigma d}{2\sigma+2} 
\|u\|_{2\sigma+2}^{2\sigma+2}$.
Applying the evolution of modified momentum Eq. \eqref{bg} and  taking expectation, we obtain 
\begin{align}\label{mod-mom}
\E\Big[e^{bt}G(u(t))\Big]&= \E\Big[ G(u_0)\Big]+2\int_{0}^{t}e^{bs}\E\Big[\tilde H(u(s)) \Big]ds-(2a-b)\int_0^t \E\Big[e^{bs}G(u(s))\Big]ds
\end{align}
To control the second term  $\E\Big[\tilde H(u(s)) \Big]$ uniformly, we take $b\le a[2-\frac {4\sigma}{\sigma d-2}]$ such that
\begin{align*}
\|\nabla u\|^2+\frac 2{b-2a}(a-\frac b{2\sigma+2})\|u\|_{2\sigma+2}^{2\sigma+2}&\ge
2H(u) +\frac {2-\sigma d}{2\sigma+2}\|u\|_{2\sigma+2}^{2\sigma+2}\\
&=\tilde H(u(s)).
\end{align*}
Then  the fact that $\tilde H(u) \le 2H(u)$ leads that
\begin{align*}
\E\Big[\frac 12e^{bt}\tilde H(u(t))\Big]+&(2a-b)\int_0^t \E \Big[\frac 12e^{bs}\tilde H(u(s))\Big]ds\\
&\le H(u_0)+\frac {1-e^{(b-2a)t}}{4a-2b}\|u_0\|^2\|f_Q\|_{L^{\infty}}.
\end{align*}
Using Gronwall inequality, we have
\begin{align*}
\int_0^t\E \Big[\frac 12e^{bs}\tilde H(u(s))\Big]ds
&\le \frac {1-e^{(b-2a)t}}{2a-b}H(u_0)
+\frac {1-e^{(b-2a)t}}{2(2a-b)^2}\|u_0\|^2\|f_Q\|_{L^{\infty}}\\
&\quad-\frac {te^{(b-2a)t}}{4a-2b} \|u_0\|^2\|f_Q\|_{L^{\infty}},
\end{align*}
which derives that
\begin{align*}
\int_0^t\Big[e^{bs}\tilde H(u(s))\Big]ds
&\le \frac 1{2a-b}\Big(2H(u_0)+\frac 1{2a-b}\|u_0\|^2\|f_Q\|_{L^{\infty}}\Big).
\end{align*}
The above estimation and 
Eq. \eqref{mod-mom} yield that
\begin{align*}
\E\Big[e^{bt}G(u(t))\Big]
&\le \E\Big[ G(u_0)\Big]-(2a-b)\int_0^t e^{bs}G(u(s))ds\\
&\quad+\frac 2{2a-b}\Big(2H(u_0)+\frac 1{2a-b}\|u_0\|^2\|f_Q\|_{L^{\infty}}\Big).
\end{align*}
Again by Gronwall inequality, 
\begin{align*}
&\E\Big[\int_0^te^{bs}G(u(s))ds\Big]\\
\le &\frac {1-e^{(b-2a)t}}{2a-b}\E\Big[G(u_0)
+\frac 2{2a-b}\Big(2H(u_0)+\frac 1{2a-b}\|u_0\|^2\|f_Q\|_{L^{\infty}}\Big)\Big].
\end{align*}
The above inequality, Eq. \eqref{bv} and the non-negativity of $V$ yield that
\begin{align*} 
&\quad \E \Big[e^{bt}V(u(t))\Big]\\
&=\E \Big[V(u_0)\Big]+4\int_0^te^{bs}\E\Big[G(u(s))\Big]ds-(2a-b)\int_{0}^\tau e^{bs}\E\Big[V(u(s))\Big]ds\\
&\le \E \Big[V(u_0)\Big]+\frac 4{2a-b}(1-e^{(b-2a)t})\E\Big[G(u_0)+\frac 2{2a-b}\Big(2H(u_0) \\
&\quad+\frac 1{2a-b}\|u_0\|^2\|f_Q\|_{L^{\infty}}\Big)\Big].
\end{align*}
Since the assumption means that $\E \Big[\frac {2a-b}4 V(u_0)+G(u_0)+\frac 2{2a-b}\Big(2H(u_0)+\frac 1{2a-b}\|u_0\|^2$ $\|f_Q\|_{L^{\infty}}\Big)\Big]< 0$, we only need to make 
\begin{align*}
\bar t \ge - \frac 1{2a-b} \ln\Bigg(\frac{\frac {2a-b}4 \E\Big[V(u_0)\Big]+\E \Big[G(u_0)+\frac 2{2a-b}\Big(2H(u_0)+\frac 1{2a-b}\|u_0\|^2\|f_Q\|_{L^{\infty}}\Big)\Big] }{\E \Big[G(u_0)+\frac 2{2a-b}\Big(2H(u_0)+\frac 1{2a-b}\|u_0\|^2\|f_Q\|_{L^{\infty}}\Big)\Big]}\Bigg).
\end{align*}
By the above inequality and the positivity of $e^{bt}$, there exists some $\bar t$ such that $\lim\limits_{t\to \bar t}\E\Big[V(u(t))\Big]=0$.
By the uncertainty principle $\|u\|^2\le \frac 2d \|\nabla u\|\big\||x|u\big\|$, we get 
\begin{align*}
\sqrt{\E\Big[\|\nabla u(t)\|^2\Big]}&\ge \frac {d \|u_0\|^2}{2e^{2a t}\sqrt{\E \Big[ \big\||x|u\big\|^2\Big]}}\\
&=\frac {d \|u_0\|^2}{2e^{(2a-b) t}\sqrt{\E \Big[ V(u(t))\Big]}}.
\end{align*}
The above estimation yields that $\sqrt{\E\Big[\|\nabla u(t)\|^2\Big]}$ goes into $\infty$ when $t$ tends to $\bar t$, which leads to a contradiction  and  finishes the proof.
\hfill\qed

\begin{rk}
The above proposition implies that when 
$$\E\Big[ V(u_0)+\frac {\sigma d-2}{a\sigma }G(u_0)+(\frac {\sigma d-2}{a\sigma})^2 H(u_0)+\frac 18(\frac {\sigma d-2}{a\sigma})^3 
\|u_0\|^2\|f_Q\|_{\infty}\Big]<0,$$
the solution of Eq. \eqref{nls} will blow up in a finite time with a positive probability.  Clearly, if the energy of $u_0$ is negative a.s., damped effect is not strong and the noise is small enough, the blow up phenomenon of the solution always happens.
\end{rk}

The blow-up condition Eq. \eqref{bcon} presents the effect of the damped term. The following result  gives a sharp time-dependent blow-up condition.  

\begin{tm}\label{blow-cor}
Let
$\sigma d>2, a\ge 0$, $u_0\in L^2(\Omega; \Sigma)\cap L^{2\sigma+2}(\Omega; L^{2\sigma+2})$,  and $\sum\limits_{k\in \N^+}\|Q^{\frac 12}e_k\|^2+\|f_Q\|_{L^{\infty}}<\infty$. 
Assume also that for some $z \ge \frac{4a\sigma}{\sigma d-2} $ and $\bar t$ such that
\begin{equation}
\begin{split}
\label{bcon1}
\E\Big[V(u_0)\Big]
&+4\bar t\E\Big[G(u_0)\Big]+(8\bar t^2+\frac 83 z\bar t^3)\E\Big[H(u_0)\Big]\\
&
+\Big(\frac 43\bar t^3+\frac 43 z\bar t^4\Big)\E\Big[\|u_0\|^2\Big]\|f_Q\|_{L^{\infty}}\le 0,
\end{split}
\end{equation}
then we have 
\begin{align*}
\PP(\tau^*(u_0)\le \bar t)>0.
\end{align*}
\end{tm}
\textbf{Proof}
The proof is similar to the proof of Proposition \ref{blow}.
Using the evolutions of the modified energy Eq. \eqref{bh}, the new energy $\tilde H$ in Proposition \ref{blow} leads that for $z\ge \frac {4a\sigma}{\sigma d-2}$, 
\begin{align*}
&\quad\int_0^t\E \Big[e^{bs}\tilde H(u(s))\Big]ds\\
&\le \frac {1-e^{-zt}}{z }\E\Big[2H(u_0)\Big]+
\frac {1-(1+zt)e^{-zt}}{z^2}\E\Big[\|u_0\|^2\Big]\|f_Q\|_{L^{\infty}}.
\end{align*}
Applying the the evolution of the modified momentum Eq. \eqref{bg}, invariance Eq. \eqref{bv}, and combining with 
the Taylor expansion of   $e^z, z\in \R$,
we get 
\begin{align*}
&\quad\E\Big[e^{bt}V(u(t))\Big]\\
&\le \E\Big[V(u_0)\Big]+\frac {1-e^{-zt}}{z}\E\Big[4G(u_0)\Big]+
\frac {1-(1+zt)e^{-zt}}{z^2}\E\Big[16H(u_0)\Big]\\
&\quad+
\frac {1-(1+zt+\frac 12(zt)^2)e^{-zt}}{z^3}
8\E\Big[\|u_0\|^2\Big]\|f_Q\|_{L^{\infty}}\\
&\le  \E\Big[V(u_0)\Big]
+4t\E\Big[G(u_0)\Big]+(8t^2+\frac 83 zt^3)\E\Big[H(u_0)\Big]\\
&\quad
+\Big(\frac 43t^3+\frac 43 zt^4\Big)\E\Big[\|u_0\|^2\Big]\|f_Q\|_{L^{\infty}}
\end{align*}
Similar arguments in Proposition \ref{blow} yield that Eq. \eqref{bcon1} is the blow-up condition.
\qed

\begin{rk}
The above blow-up condition is sharp in the sense that 
if $z\to 0$, the above condition can be degenerated to the blow-up condition in  the conservative case, i.e. for some $\bar t$
\begin{align*}
\E\Big[V(u_0)\Big]+4\E\Big[G(u_0)\Big]\bar t+8\E\Big[H(u_0)\Big]\bar t^2+\frac 43\bar t^3\|f_Q\|_{L^{\infty}}\E\Big[M(u_0)\Big]<0.
\end{align*}
 Since the blow-up condition for Eq.\eqref{nls} in supercritical case is similar
to the condition of Theorem 4.1 in \cite{DD05}, it is possible to apply the 
skills and arguments in \cite{DD05} to get a stronger result that  $\PP(\tau^*(u_0,a)>t)>0$, for any $t$ and $u_0$, $u_0\ne 0$ under some assumptions on $u_0$, $d$ and $Q^{\frac 12}$.
\end{rk}

\begin{rk}
The blow-up condition in critical case can not be obtained by the method in Proposition \ref{blow},
we only get some necessary condition in Remark \ref{nes}. It is our future work to study the blow-up phenomena of the solution and show the sufficient blow-up condition  for damped stochastic NLS equations in critical case.
\end{rk}

It seems  that when $a$ becomes larger, the blow-up time becomes longer  and that  when $a$ goes to $\infty$, the blow-up condition is not satisfied. Indeed, we can show that 
when the damped effect is large enough, the damped effect can prevent the  blow-up of the solution with high probability. The key of the proof is using the  infinite dimensional Doss--Sussman type transformation in \cite{BRZ16, BRZ17}.
In the  following theorem, 
we assume that the noise satisfies $\sum_{k\in \N^+}\| Q^{\frac 12}e_k\|_{C_b^2}<\infty$ and 
the following decay condition
\begin{align*}
\lim_{x\to \infty}\eta(x)(|Q^{\frac 12}e_m(x)|+|\nabla Q^{\frac 12} e_m(x)|+|\Delta Q^{\frac 12}e_m(x)|)=0,
\end{align*}
where $\eta(x)=1+|x|^2$ if  $d\ne 2$, and 
$\eta(x)=(1+|x|^2)(\ln(2+|x|^2))^2$ if  $d= 2$.
Under these assumptions, the local well-posedness is obtained  in \cite{BRZ16}. 
We also remark when $\{e_k\}_{k\in N^+}$ is an orthonormal basis 
of $\HH$, the decay condition  natural holds. 
In this case,  $e_k$ can be chosen as the $k$-th Hermite function in $\HH$, and meanwhile as the $k$-th eigenvector of the operator $Q^{\frac 12}.$

\begin{tm}\label{noblow}
Assume that $a>0$, $\sigma d \ge 2$ and 
$\sum_{k}\|\Delta Q^{\frac 12} e_k\|_{L^{\infty}}^2<\infty$. Then for any $u_0\in \HH^1$ and 
$0<T<\infty$, we have
\begin{align*}
\lim_{a\to \infty}\PP(\tau^*(u_0,a)> T)=1.
\end{align*}
\end{tm}
\textbf{Proof}
We apply the rescaling transformation $v(t)=e^{at-\bi W(t)}u(t)$ to Eq. \eqref{nls}
and  get the following random  partial differential equation
\begin{align}\label{rnls}
dv&=\bi\exp\Big(at -\bi W(t)\Big)\Delta u(t)dt+\bi \exp\Big(at -\bi W(t)\Big)|u|^{2\sigma}udt\\\nonumber
&=\bi\Big( \Delta +2\bi\nabla(W(t))\cdot \nabla +|\nabla W(t)|^2+\bi \Delta W(t)\Big)v \\\nonumber
&\quad+\bi\exp\Big(-2a\sigma t\Big)|v|^{2\sigma}vdt\\\nonumber
&:= A(t)vdt
+\bi\exp\Big(-2a\sigma t\Big)|v|^{2\sigma}vdt.
\end{align}
By Lemma 2.4 in \cite{BRZ16}, the solutions of Eq. \eqref{rnls} and Eq. \eqref{nls}
are equivalent. 
Let $r=  \frac {4\sigma+4}{\sigma d}$ such that $(2\sigma+2, r)$ is a Strichartz pair. Next we recall the proof 
of  the local well-posedness and set
\begin{align*}
\mathcal X_R^\tau:=\Big\{ v\in C(0,\tau; \HH^1)\cap L^{r}(0,\tau; W^{1,2\sigma+2})\Big| \|v\|_{C(0,\tau; \HH^1)}+\|v\|_{L^{r}(0,\tau; W^{1,2\sigma+2})}\le R\Big\}. 
\end{align*}
Considering the solution map $G$ of Eq. \eqref{rnls},
by the random Strichartz estimate and similar arguments in \cite{BRZ16}, 
we obtain
\begin{align*}
&\|G(v)\|_{C(0,\tau; \HH^1)}
+\|G(v)\|_{L^{r}(0,\tau; W^{1,2\sigma+2})}\\
&\qquad\le 
2C_\tau \|u_0\|_{\HH^1}+2R^{2\sigma+1}C_{S}C_\tau \|\exp(-2a\sigma t)\|_{L^{q}(0,\tau)},\\ 
&\|G(v)-G(w)\|_{C(0,\tau; \HH)}
+\|G(v)-G(w)\|_{L^{r}(0,\tau; L^{2\sigma+2})}\\
&\qquad \le 4C_{S}C_\tau R^{2\sigma}\|\exp(-2a\sigma t)\|_{L^{q}(0,\tau)} \|v-w\|_{L^{r}(0,\tau; L^{2\sigma+2})},
\end{align*} 
where $v,w \in \mathcal X_R^\tau$, $C_{S}=(2\sigma+1)D^{2\sigma}$, $D$ is the Sobolev embedding coefficient form $L^{2\sigma+2}$ to
$\HH^1$, $C_\tau$ is the random Strichartz estimate coefficient, $q>1$ and $\frac 1q=1-\frac 2 r>0$.
Now we take $R=4C_\tau \|u_0\|_{\HH^1}$ and $$\tau(a) =\inf\Big\{t>0\Big|\frac {2C_{S}C_t}{a\sigma}R^{2\sigma}>1\Big\}$$
such that  the mapping $G:X_R^{\tau}\longrightarrow X_R^{\tau}$ 
has a fixed point in the Banach space $\Big(X_R^{\tau}, \|\cdot\|_{C(0,\tau; \HH)}+\|\cdot \|_{{L^{r}(0,\tau; L^{2\sigma+2})}}\Big)$, which implies the local well-posedness of Eq. \eqref{rnls}.

Now,  we aim to show that $\lim_{a\to \infty} \PP(\tau^*(u_0,a)>T)=1$. 
Based on the result in \cite[Lemma 2.7]{BRZ16} that $C_t$, $t>0$ is $\FFF_t$ measurable, increasing
and continuous,    
the definition of $\tau^*(u_0, a)$
and the a.s. boundedness of $C_t$ yield that
\begin{align*} 
\lim_{a\to \infty} \PP(\tau^*(u_0,a)> T)
&\ge 
\lim_{a\to \infty} \PP(\tau(a)> T)\\
&=\lim_{a\to \infty} \PP( \frac {2C_{S}C_t}{a\sigma}R^{2\sigma}\le 1,  \, t \in [0,T])\\
&\ge \lim_{a\to \infty} \PP( \frac {2C_{S}C_T}{a\sigma}R^{2\sigma}\le 1)\\
&\ge \lim_{a\to \infty} \PP( \frac {2C_{S}C_T}{\sigma}R^{2\sigma}\le a)
=1.
\end{align*}

\qed

\begin{rk}
If the noise is space-independent or disappears,
applying the arguments in Theorem \ref{well} and Theorem \ref{noblow}, one can get global existence 
and blow-up results of the solutions for the damped stochastic NLS equation. 
\end{rk}

\section{Conclusions}\label{sec;5}

In this paper, we consider the influence of both damped term and noise on the stochastic nonlinear Schr\"odinger equation driven by multiplicative noise. 
We first show the global existence of the unique solution for damped stochastic NLS equation in critical case and study the exponential integrability  and the continuous dependence on the initial data of the solution. 
Then based on the modified variance identity, we deduce a sharp  blow-up condition for damped stochastic NLS  equation in supercritical case.
Moreover, we prove that the large damped effect can prevent 
the blow-up with high probability.

\bibliography{bib}
\bibliographystyle{plain}

\end{document}